\documentclass[12pt]{amsart}
\usepackage{amsmath,amssymb,amsthm,amsxtra}
\numberwithin{equation}{section}

\newtheorem{thm}{Theorem}[section]
\newtheorem{cor}[thm]{Corollary}

\newtheorem{lem}[thm]{Lemma}

\newtheorem{prp}[thm]{Proposition}
\newtheorem{rmk}[]{Remark}%[section]

\newcommand{\im}{\operatorname{Im}}
\newcommand{\re}{\operatorname{Re}}

\newcommand{\mes}{\operatorname{mes}}
\def\sgn{\operatorname{sgn}}
\newcommand{\rmd}{\mathrm{d}}
\newcommand{\rme}{\mathrm{e}}
\newcommand{\rmi}{\mathrm{i}}
\newcommand{\loc}{\mathrm{loc}}

\newcommand{\wto}{\to\kern-1.22em\raise.9ex
                  \hbox{\small\rm w}\hskip.7em} % weakly tends to
\newcommand{\wwto}{\to\kern-1.2em\raise.9ex
                  \hbox{\small\rm w}\hskip.7em} % weakly tends to
\def\lll{\lambda}

\newcommand{\calD}{\mathcal{D}}

\newcommand{\mbC}{\mathbb{C}}

\newcommand{\mbN}{\mathbb{N}}

\newcommand{\mbR}{\mathbb{R}}

\newcommand{\aaa}{\alpha}
\newcommand{\bbb}{\beta}

\begin{document}

%\addtolength{\textwidth}{2in}
%\addtolength{\textheight}{3in}

\title[Non-real eigenvalues of indefinite Sturm-Liouville problems]
{A priori bounds and existence of non-real eigenvalues of
indefinite Sturm-Liouville problems}

\author[J.\,Qi and S.\,Chen]
       {Jiangang Qi and Shaozhu Chen}
\address{Department of Mathematics, Shandong University(Weihai)\\
        Weihai 264209, P.R. China}
\email{qjg816@163.com; szchen@sdu.edu.cn}

\begin{thanks} {This research was partially supported by the NSF of China (Grant 11271229),
the NSF of Shandong Province (Grants ZR2012AM002 and ZR2011AQ002) and the Teaching Research
and Teaching Reform Project of Shandong University, Weihai (Grant A201005). }
\end{thanks}

%\date{\today}

\begin{abstract} The present paper gives a priori bounds
on the possible non-real eigenvalues of
regular indefinite Sturm-Liouville problems
and obtains sufficient conditions for such problems
to admit non-real eigenvalues.
\medskip

\vspace{04pt}

\end{abstract}

\maketitle

%S1
%%%%%%%%%%%%%%%%%%%%%%%%
\section{Introduction}
\setcounter{equation}{0}
%%%%%%%%%%%%%%%%%%%%%%%%

The present paper is concerned with the indefinite  spectral problem
\begin{equation}                                         \label{prob}
-y''+q y=\lll w y, \ y(-1)=y(1)=0, \; \text{ in }\; L^2_{|w|}[-1,1]
\end{equation}
under the standing hypothesis that $q$ and $w$ are real-valued functions
satisfying
\begin{equation}                                         \label{cd}
w(x)\ne 0\ \text{a.e. on } [-1,1], \; q, w\in L^1[-1,1],
\end{equation}
and $w(x)$ changes sign on $[-1,1]$.
The indefinite problem \eqref{prob} has discrete, real eigenvalues,
unbounded from both below and above,
and may also admit non-real eigenvalues.
Such problems occur in certain physical models, particularly in
transport theory and statistical physics.
The indefinite nature of the problem was
noticed by Haupt \cite{Hau} and Richardson
\cite{rich1912} at the beginning
of the last century.
For a review of the early work in this direction, see \cite{min2}.

As a simple example of \eqref{prob}, the Richardson problem  \cite{rich1918}
\begin{equation}\label{son}
-y''-\mu y=\lll {\rm sgn}(x)y,\ x\in[-1,1],\  y(-1)=0=y(1)
\end{equation}
was studied by many authors, such as
Turyn \cite{tur}, Atkinson and Jabon \cite {atkinson},
Fleckinger and Mingarelli \cite{fle}, and
P.\,Binding and H.\,Volkmer \cite{Bind}.
For the indefinite problem \eqref{prob},
non-real eigenvalues might appear only if
the corresponding right-definite
problem
\begin{equation}                                         \label{15}
-y''+q y=\lll |w| y,\; y(-1)=y(1)=0\;\text{ in }\; L^2_{|w|}[-1,1]
\end{equation}
has negative eigenvalues, namely, here holds the following
result.

\begin{prp}\label{p1} {\rm(cf. \cite[Theorem 2, p.\,523]{min1} and
\cite[Corollary 1.7]{Langer})}
If the problem \eqref{15} has $n$ negative eigenvalues,
then the problem \eqref{prob} has at most $2n$ non-real eigenvalues.
\end{prp}

Although the upper bound given in Proposition \ref{p1} is sharp
\cite{rich1912, Trunk},
determining a priori bounds and the exact number of non-real eigenvalues
are still  difficult and  interesting open problems in
Sturm-Liouville theory (see \cite{min2} and \cite[p.\,126]{Zet}).
Recently, by means of the operator theory in Krein spaces,
Behrndt, Katatbeh and Trunk
\cite[Theorem 2.3, Corollary 2.4]{Trunk} have given
sufficient conditions for the existence of non-real eigenvalues
of the singular indefinite Sturm-Liouville operator
\begin{equation}\label{ph}
(Af)(x):= \sgn(x)(-f''(x) + V (x)f(x))=\lll f(x), \ x\in\mbR,
\end{equation}
%An approach used in \cite{Trunk} is essentially based on the fact that the
%essential spectrum of the singular operator  $A$ lies on the whole real axis,
%and hence, cannot be applied to the regular cases.
and if $V\in L^\infty(\mbR)$, Behrndt, Philipp and Trunk
\cite[Theorem 4,2]{Trunk3} have obtained explicit bounds
on the non-real eigenvalues of \eqref{ph} in terms of $V$.

\medskip
In the present paper, we will first obtain a priori bounds for possible
non-real eigenvalues and then find sufficient conditions for the existence of
non-real eigenvalues of \eqref{prob}.
These results will answer or partially answer several open problems
posed in \cite{min2}.
We state these results in this section and prove them in
Sections 2 and 3.

Denote by $\|\cdot\|_p$ the norm of the space $L^p[-1,1]$ and by
$\|\cdot\|_C$ the maximum norm of $C[-1,1]$.
If $xw(x)>0$ a.e. on $[-1,1]$, we set
\begin{equation}\label{m1}
S_1(\varepsilon)=\{x\in[-1,1]:\ xw(x)<\varepsilon\},\
m_1(\varepsilon)=\mes S_1(\varepsilon).
\end{equation}
If $w\in AC_{\loc}[-1,1]$, $w'\in L^2[-1,1]$, we set
\begin{equation}\label{m2}
S_2(\varepsilon)=\{x\in[-1,1]:\ w^2(x)<\varepsilon\},\
m_2(\varepsilon)=\mes S_2(\varepsilon).
\end{equation}
A value of $x$ about which $w(x)$ changes its sign will be called a
{\it turning point} \cite{min1}.
If $w(x)$ has only one turning point, we will obtain the following
a priori bounds for possible non-real eigenvalues.

\begin{thm}\label{thm0}
Suppose that $\lambda$ is, if it exists, a non-real eigenvalue of \eqref{prob}.
If $xw(x)>0$ a.e. on $[-1,1]$, then
\begin{equation}\label{bd1.mu=0}
|\re\lll|\le \frac{4}{\varepsilon_1}\big(\|q_-\|_1 +4\|q_-\|_1^2\big),\
|\im\lll|\le \frac4{\varepsilon_1}\|q_-\|_1.
\end{equation}
where $\varepsilon_1>0$ satisfies $8\|q_-\|_1^2m_1(\varepsilon_1)<1$ and $q_-(x)=-\min\{0,q(x)\}$.
\end{thm}
In the case where $w(x)$ is allowed to have more
turning points, we will obtain

\begin{thm}\label{thm0'}
Suppose that $\lambda$ is, if it exists, a non-real eigenvalue of \eqref{prob}.
If $w\in AC[-1,1]$ and $w'\in L^2[-1,1]$, then
\begin{equation}\label{bd2.mu=0}
|\re\lll|\le \frac{8}{\varepsilon_2}\|q_-\|_1^2\left(3\|w\|_C+\|w'\|_2\right),\
|\im\lll|\le \frac8{\varepsilon_2}\|w'\|_2\|q_-\|^2_1.
\end{equation}
where $\varepsilon_2>0$ is chosen such that
$8\|q_-\|_1^2m_2(\varepsilon_2)<1.$
\end{thm}

In the particular case where $q\ge 0$, we see by Theorems \ref{thm0}
and \ref{thm0'}
that \eqref{prob} has no any non-real eigenvalues, which is in accordance
with the conclusion in Proposition \ref{p1} since now \eqref{15}
does not have any negative eigenvalues.

In what follows, we impose the symmetry conditions on $q$ and $w$, namely,
\begin{equation}\label{cdsym}
q(x)=q(-x), \ w(-x)=-w(x).
\end{equation}
In this case, more accurate a priori bounds on imaginary eigenvalues can be found
if $q$ is bounded below and $w$  keeps away from zero.

\begin{thm}\label{cor1}
Suppose that \eqref{cdsym} holds and $xw(x)>0$ a.e. on $[-1,1]$.
If, for some $q_0<0$ and $w_0>0$,
\begin{equation}                                 \label{bddqw}
q(x)\ge q_0, \ |w(x)|\ge w_0\ \text{ a.e. }\;x\in[-1,1],
\end{equation}
then for any possible pure imaginary eigenvalue  $\lll$ of \eqref{prob},
there holds
\begin{equation}\label{in4mu=0}
|\im\lll|\le \frac{4(-q_0)^{3/2}}{w_0}.
\end{equation}
\end{thm}

In view of \eqref{cdsym}, using the spectral theory of operators
in Krein spaces, we obtain
an existence result for non-real eigenvalues
of the indefinite  problem \eqref{prob}.

\begin{thm}\label{thm1}
Let \eqref{cdsym} be fulfilled.
If the eigenvalue problem
\begin{equation}\label{w=1}
-y''+q(x)y=\lll y, \  y(-1)=y(1)=0
\end{equation}
has one negative eigenvalue and the rest
eigenvalues are all positive,
then \eqref{prob} has exactly two purely
imaginary eigenvalues.
\end{thm}

Immediate consequences of Proposition \ref{p1}, Theorems \ref{cor1} and \ref{thm1}
are the existence and bounds for non-real eigenvalues of
Richardson problem.

\begin{cor}\label{thm12}
For $\mu\in\left(\frac{\pi^2}{4}, \pi^2\right)$,
the Richardson eigenvalue problem \eqref{son}
has  exactly  two purely imaginary  eigenvalues whose moduli are bounded by $4\mu^{3/2}$.
\end{cor}

\begin{rmk}  Theorems \ref{thm0}--\ref{thm1} can be generalized
to the problem
$$\left\{\begin{aligned}
&-(p(x)y')'+q(x)y=\lambda w(x)y,\\
& \alpha_1 y(-1)+\bbb_1 y'(-1)=0,\\
& \alpha_2 y(1)+\bbb_2y'(1)=0,
\end{aligned}\right.
$$
where $p(x)>0$ a.e. on $[-1,1]$, $1/p\in L^1[-1,1]$, $\alpha_j,\beta_j\in\mbR$
for $j=1,2$ and $\alpha_1 \bbb_2+\alpha_2\bbb_1=0$,
but we do not pursue this here.
\end{rmk}

%Following this section, Theorems \ref{thm0}, \ref{thm0'} and \ref{cor1} will be proved in Section 2
%and Theorem \ref{thm1} will proved in Section 3.

%%%%%%%%%%%%%%%%%%%%%%%%%%%%%%%%%%%%%%%%%%%%%%%%%
\section{A priori bounds of non-real eigenvalues}
%%%%%%%%%%%%%%%%%%%%%%%%%%%%%%%%%%%%%%%%%%%%%%%%%

In this section we will prove Theorems \ref{thm0}, \ref{thm0'} and \ref{cor1}.
\medskip

\noindent{\bf The proof of Theorem \ref{thm0}.}
Let $\lambda$ be a non-real eigenvalue of \eqref{prob} and $\phi(x)$
the corresponding eigenfunction  with  $\|\phi\|_2=1$.
Multiplying both sides of $-\phi''+q\phi=\lll w\phi$
by $\overline{\phi}$
and integrating over the interval $[x,1]$ we have
\begin{equation}\label{int.xto1}
(\phi'\overline{\phi})(x)+\int^1_x|\phi'|^2+\int^1_x q|\phi|^2=\lll\int^1_x w|\phi|^2.
\end{equation}
Separating the real and imaginary parts of  both sides of \eqref{int.xto1}
yields
\begin{eqnarray}
\re\lll\int^1_x w|\phi|^2
&=&\re(\phi'\overline{\phi})(x)+\int^1_x|\phi'|^2+\int^1_x q|\phi|^2,\label{312}\\
\im\lll\int^1_x w|\phi|^2&=&\im(\phi'\overline{\phi})(x). \label{313}
\end{eqnarray}
We will use \eqref{312} and \eqref{313} to estimate $\re\lambda$ and $\im\lambda$.
To do this, let $x=-1$ in \eqref{313}. From $\im\lll\not=0$ and $\phi(-1)=0$,
we have $\int^1_{-1} w|\phi|^2=0$,
and hence, by \eqref{312},
\begin{equation}\label{3140}
\int_{-1}^1(|\phi'|^2+q|\phi|^2)=0.
\end{equation}
Set $Q(x)=\int^x_{-1} q_-(t)\rmd t$. Then $\max |Q(x)|\le \|q_-\|_1$ and
$$
\begin{aligned}
\int^1_{-1}q_-(x)|\phi(x)|^2\rmd x
&=\int^1_{-1} Q'(x)|\phi(x)|^2\rmd x\\
&-2\re\left(\int^1_{-1} Q(x)\phi'(x)\overline{\phi(x)}\rmd x\right),
\end{aligned}
$$
which, together with $\|\phi\|_2=1$, yields that
\begin{equation}\label{in1}
\int^1_{-1}q_-|\phi|^2\le 2\|q_-\|_1 \int^1_{-1} |\phi'||\phi|
\le 2\|q_-\|_1\|\phi'\|_2\le 2\|q_-\|^2_1+\frac{1}{2}\|\phi'\|^2_2.
\end{equation}
Then, from \eqref{3140}, we get
%$$
%\|\phi'\|_2^2\le\int_{-1}^1q_-|\phi|^2\le2\|q_-\|_1^2+\frac12\|\phi'\|_2^2,
%$$
%and hence,
\begin{equation}\label{316}
\|\phi'\|_2^2\le4\|q_-\|^2_1,\;\;
\int_{-1}^1q_-|\phi|^2\le4\|q_-\|^2_1.
\end{equation}
%Therefore, using \eqref{in1} again we have
%\begin{equation}\label{injia}
%\left|\int^1_{-1}q_-|\phi|^2\right|\le 2\|q_-\|_1\|\phi'\|_2\le 4\|q_-\|^2_1.
%\end{equation}
From $\phi(x)=\int_{-1}^x\phi'(t)\rmd t$,
by Cauchy-Schwarz inequality, we have
$$
|\phi(x)|^2=\left|\int_{-1}^x\phi'(t)\rmd t\right|^2
\le(x+1)\int_{-1}^x|\phi'(t)|^2\rmd t\le\int_{-1}^0|\phi'|^2
\le\|\phi'\|_2^2
$$
for $-1\le x\le0$.
From $\phi(x)=-\int_x^1\phi'(t)\rmd t$, one similarly proves
$|\phi(x)|^2\le\|\phi'\|_2^2$ for $x\in[0,1]$, and so,
\begin{equation}\label{317}
|\phi(x)|^2%=\left|\int^x_{-1}\phi'(t)\rmd t\right|^2
\le\|\phi'\|_2^2, \;\; x\in[-1,1].
%\le 4|\mu|+8\|q_-\|_1^2
\end{equation}

Since $xw(x)>0$, a.e. on $[-1,1]$, one can find $\varepsilon_1>0$
such that $8\|q_-\|_1^2m_1(\varepsilon_1)<1$,
where $m_1(\varepsilon)$ is defined in \eqref{m1}.
Using $\int_{-1}^1w|\phi|^2=0$, from \eqref{316} and \eqref{317}, we have
\begin{equation}\label{in2}
\aligned
\int_{-1}^1\int_x^1w(t)|\phi(t)|^2\rmd t\,\rmd x
%=\int_{-1}^1 w(t)|\phi(t)|^2\int^t_{-1}\rmd x\,\rmd t
&%=\int_{-1}^1(t+1)w(t)|\phi(t)|^2\rmd t
=\int^1_{-1} xw(x)|\phi(x)|^2\rmd x\\
&\ge\varepsilon_1\left(\int^1_{-1}|\phi(x)|^2\rmd x
  -\int_{S_1(\varepsilon_1)}|\phi(x)|^2\rmd x\right)\\
&\ge\varepsilon_1\left[1-4\|q_-\|_1^2 m_1(\varepsilon_1)\right]
 \ge\frac{\varepsilon_1}{2}.
\endaligned
\end{equation}

Set $q_+(x)=\max\{0,q(x)\}$. Then $q=q_+-q_-$ and $|q|=q_++q_-=q+2q_-$.
Repeatedly using \eqref{3140}, we have
$$
\aligned
\left|\int_{-1}^1\int_x^1(|\phi'|^2+q|\phi|^2)\rmd t\rmd x\right|
&%=\left|\int_{-1}^1(x+1)(|\phi'|^2+q|\phi|^2)\rmd x\right|
 =\left|\int_{-1}^1x(|\phi'|^2+q|\phi|^2)\rmd x\right|\\
&\le\int_{-1}^1(|\phi'|^2+q|\phi|^2+2q_-|\phi|^2)\rmd x\\
&=2\int_{-1}^1q_-|\phi|^2\rmd x.
\endaligned
$$
Now, by \eqref{316}, the integration of \eqref{312} gives
$$
\aligned
|\re\lll|\int^1_{-1}\int^1_x w|\phi|^2
&=\left|\int_{-1}^1\re(\phi'\overline{\phi})\rmd x
  +\int_{-1}^1\int_x^1(|\phi'|^2+q|\phi|^2)\rmd t\rmd x\right|\\
%&\le \int^1_{-1}|\phi'||\phi|+\int^1_{-1}\left(|\phi'|^2+|q||\phi|^2\rmd x\right)\\
%&\le \|\phi'\|_2+\int^1_{-1}\left(|\phi'|^2+(q+2q_-)|\phi|^2\right)\\
&\le\|\phi'\|_2+2\int_{-1}^1 q_-|\phi|^2\rmd x\le 2\|q_-\|_1+8\|q_-\|_1^2.
\endaligned
$$
Therefore, in view of \eqref{in2}, we conclude that
\begin{equation}\label{real}
|\re\lll|\le \frac4{\varepsilon_1}\big(\|q_-\|_1 +4\|q_-\|_1^2\big).
\end{equation}
Moreover, integrating \eqref{313} and
using \eqref{in2} and \eqref{316}, we have
\begin{equation}
\frac{\varepsilon_1}2|\im\lambda|
\le|\im\lambda|\int_{-1}^1\int_x^1w|\phi|^2
=\left|\int_{-1}^1\im(\phi'\overline{\phi})\right|
%\le \int_{-1}^1|\phi'\overline{\phi}|
\le\|\phi'\|_2\le2\|q_-\|_1,
\end{equation}
and \eqref{bd1.mu=0} follows immediately. %one sees that $|\im\lll|\le 4\|q_-\|_1/{\varepsilon_1}.$
This  completes the proof of Theorem \ref{thm0}.\qed

\medskip

\noindent {\bf The Proof of Theorem \ref{thm0'}.}
Let $\lambda$ be a non-real eigenvalue of \eqref{prob} and $\phi$
the corresponding eigenfunction  with  $\|\phi\|_2=1$.
In this case we still can make use of \eqref{int.xto1}, \eqref{312} and \eqref{313}.
From \eqref{313}, since $\im\lambda\ne0$, one sees that
$\int^1_{-1}w |\phi(x)|^2\rmd x=0$.
Thus, \eqref{3140}, \eqref{316} and \eqref{317}
hold, and particularly,
\begin{equation}\label{3170}
|\phi(x)|^2\le\|\phi'\|_2^2,\ x\in[-1,1],\ \|\phi'\|_2^2\le\int^1_{-1}q_-|\phi|^2\le 4\|q_-\|_1^2.
\end{equation}
Multiplying $-\phi''+q\phi=\lll w\phi$ by $w\overline{\phi}$
and integrating by parts, we get
\begin{equation}\label{eq}
\int^1_{-1}w|\phi'|^2+\int^1_{-1}w'\phi'\overline{\phi}
+\int^1_{-1}w q|\phi|^2=\lll\int^1_{-1} w^2|\phi|^2.
\end{equation}
Separating the real and imaginary parts of the both sides of \eqref{eq}
yields
\begin{eqnarray}
\qquad\re\lll\int^1_{-1} w^2|\phi|^2
&=&\re\left(\int^1_{-1}w'\phi'\overline{\phi}\right)
+\int^1_{-1}w(|\phi'|^2+q|\phi|^2),\label{wc1}\\
\im\lll\int^1_{-1} w^2|\phi|^2
&=&\im\left(\int^1_{-1}w'\phi'\overline{\phi}\right). \label{wc2}
\end{eqnarray}
Now, using \eqref{3170}, $|q|=q+2q_-$ and $\int^1_{-1}q|\phi|^2=-\int^1_{-1}|\phi'|^2<0$,
we obtain
\begin{equation}\label{eq01}
\begin{aligned}
\left|\int^1_{-1}w|\phi'|^2\right|
&\le \|w\|_C\|\phi'\|_2^2\le4\|w\|_C\|q_-\|^2_1,\\
\left|\int^1_{-1}w q|\phi|^2\right|&\le\|w\|_C\int^1_{-1}|q||\phi|^2
\le 8\|w\|_C\|q_-\|^2_1,\\
\left|\int^1_{-1}w'\phi'\overline{\phi}\right|
&\le\|\phi'\|_2\|w'\|_2\|\phi'\|_2
\le 4\|w'\|_2\|q_-\|^2_1.
\end{aligned}
\end{equation}

Recall that $m_2(\varepsilon_2)=\mes S_2(\varepsilon_2)$ defined in \eqref{m2}
and $w^2(x)\ge\varepsilon_2$ on the set
$\Omega(\varepsilon_2):=[-1,1]\setminus S_2(\varepsilon_2)$.
Then $8\|q_-\|^2_1 m(\varepsilon_2)<1$ yields that
\begin{equation}\label{inn2}
\aligned
\int^1_{-1} w^2(x)|\phi(x)|^2\rmd x
&\ge \varepsilon_2 \int_{\Omega(\varepsilon_2)}|\phi|^2
=\varepsilon_2\left(1-\int_{S(\varepsilon_2)} |\phi|^2\right)\\
&\ge\varepsilon_2\left(1-4\|q_-\|^2_1 m(\varepsilon_2)\right)
\ge\frac{\varepsilon_2}{2},
\endaligned
\end{equation}
which, together with \eqref{wc1}, \eqref{wc2} and \eqref{eq01}, gives
\eqref{bd2.mu=0} and completes the proof. \qed
\medskip

Under the conditions  \eqref{cd} and \eqref{cdsym},
it is easy to see that if $\lambda\in\mbC$ is an
eigenvalue of \eqref{prob} with an eigenfunction $\phi$,
then $-\overline{\lambda}$ is an eigenvalue of \eqref{prob} with
the eigenfunction $\overline{\phi(-\cdot)}$.
Thus, if $\lll=\rmi\aaa$ with $\alpha\in\mbR$, then
$\overline{\phi(-x)}=C\phi(x)$ for some $C\not=0$
since the geometric multiplicity is one.
Then it follows that $|C|=1$ from $\overline{\phi(0)}=C\phi(0)$,
$\overline{\phi'(0)}=-C\phi'(0)$,
and $|\phi(0)|+|\phi'(0)|\ne0$.
To sum up, we have
\begin{lem}\label{lem21}
Let \eqref{cd} and \eqref{cdsym} hold.
If $\lambda\in\mbC$ is an eigenvalue of \eqref{prob} with
an eigenfunction $\phi$, then $-\overline{\lambda}$
is an eigenvalue of \eqref{prob} with the eigenfunction
$\overline{\phi(-\cdot)}$.
Particularly, if $\lll=\rmi\aaa$ with $\alpha\in\mbR$ and
$\alpha\ne0$, then $\overline{\phi(-\cdot)}=C\phi$ for some
$C\in\mathbb C$ with $|C|=1$.
\end{lem}

\noindent {\bf The Proof of Theorem \ref{cor1}.} Let $\phi$ be an
eigenfunction corresponding to $\lambda=\rmi\alpha$ with $\|\phi\|_2=1$.
It follows from Lemma \ref{lem21} that
there exists an $\omega\in[0,2\pi)$ such that
$\overline{\phi(-x)}=\rme^{\rmi\omega}\phi(x)$
and $-\overline{\phi'(-x)}=\rme^{\rmi\omega}\phi'(x)$.
So, $|\phi(x)|$ and $|\phi'(x)|$ are even functions.
We see that \eqref{int.xto1}-\eqref{3140}  %, \eqref{312} and \eqref{313}
hold for this $\phi$.
Similarly to \eqref{317}, we have
\begin{equation}\label{phix}
|\phi(x)|^2\le(x+1)\int_{-1}^x|\phi'(t)|^2\rmd t\le\int_{-1}^0|\phi'(t)|^2\rmd t
=\frac12\|\phi'\|_2^2,\;x\in[-1,0],
\end{equation}
since $|\phi'(x)|$ is even.
Actually, \eqref{phix} is true for $x\in[-1,1]$
since $|\phi(x)|$ is even.

Since $q(x)\ge q_0$ on $[-1,1]$, it follows from \eqref{3140} and $\|\phi\|_2=1$ that
$\|\phi'\|_2^2=-\int_{-1}^1q|\phi|^2\le-q_0$,
and then the integration of \eqref{313} produces
\begin{equation}
\label{3190}
|\im\lambda|\Big|\int_{-1}^1\int_x^1w|\phi|^2\Big|
=\Big|\int_{-1}^1\im(\phi'\overline{\phi})\Big|
\le  %\|\phi'\|_2\|\phi\|_2=
\|\phi'\|_2\le(-q_0)^{1/2}.
\end{equation}
Let $\delta=1/(-2q_0)$.
By \eqref{phix}, we have
$1=\int_{-1}^1|\phi|^2\le\|\phi'\|_2^2\le-q_0$ and
\begin{equation}\label{31900}
\begin{aligned}
\Big|\int_{-1}^1\int_x^1w|\phi|^2\Big|
&=\int^1_{-1} xw(x)|\phi|^2\rmd x
\ge w_0\int^1_{-1} |x||\phi|^2\rmd x\\
&\ge w_0\delta\int_{|x|\ge\delta}|\phi|^2
= w_0\delta\left(1-\int^\delta_{-\delta}|\phi|^2\right)\\
&\ge w_0\delta(1-\delta(-q_0))
=-\frac{w_0}{4q_0}.
\end{aligned}
\end{equation}
Now, \eqref{in4mu=0} follows from \eqref{3190} and \eqref{31900}.
The proof is complete. \qed

%%%%%%%%%%%%%%%%%%%%%%%%%%%%%%%%%%%%%%%%%%%%%%
\section{Existence of non-real eigenvalues}
%%%%%%%%%%%%%%%%%%%%%%%%%%%%%%%%%%%%%%%%%%%%%%
In this section we prove Theorem \ref{thm1} and
in the proof we will use the following result
which was proved, e.g., in \cite{Langer} and \cite{bind-moller}.

\begin{lem}\label{lem20}{\rm(}cf. \cite[Proposition 2.6]{Langer}{\rm)}
If $w_j\in L^1[-1,1]$ and $w_j(x)>0$ a.e. on $[-1,1]$
for $j=1,2$, then the two eigenvalue problems
 \begin{equation}                                          \label{wj}
-y''+q(x)y=\lll w_j(x) y, \  \ y(-1)=y(1)=0,\;\;j=1,2
\end{equation}
%problems in \eqref{wj} %$H_1$ and $H_2$
have the same number of negative eigenvalues.
\end{lem}

%In the following proof we will use the natation of negative square of a self-adjoint operator
%in Krein spaces.

Let $K$ be the Krein space $L^2_{|w|}[-1,1]$, equipped with the indefinite inner product
\begin{equation}\label{inner}
[f,g]=\int^1_{-1}f(x)\overline{g(x)}w(x)\rmd x,\ f,g\in L^2_{|w|}[-1,1]
\end{equation}
and $T$ a self-adjoint operator in $K$ with domain $\calD(T)$
\cite{Trunk2,Trunk, Langer}.
We say that the operator $T$ has $k$ negative squares, $k\in\mbN_0$,
if there exists a $k$-dimensional subspace $X$ of $K$ in
$\calD(T)$ such that $[Tf,f]<0$ if $f\in X$ and $f\not=0$, but no $(k+1)$-dimensional subspace
with this property.

\medskip

\noindent {\bf The Proof of Theorem \ref{thm1}.} Let $A$ and $B$
be the operators associated with $-y''+q(x)y=\lll w(x) y$ and
$-y''+q(x)y=\lll |w(x)|y$ with the Dirichlet boundary conditions,
respectively.
Then $B$ is self-adjoint with respect to the definite inner product
$$
(f,g)=\int^1_{-1}f(x)\overline{g(x)}|w(x)|\rmd x,\ f,g\in L^2_{|w|}[-1,1]
$$
and $A$ is self-adjoint with respect to the
indefinite inner product \eqref{inner}.

It follows from Lemma \ref{lem20} and the assumption in Theorem \ref{thm1} that
$B$ has one negative
eigenvalue and the rest are positive, and hence, $A$ has  exactly one
negative square since $[Af,f]=(Bf,f)$ and $0$ is a resolvent point of $A$.
It is well known (see, e.g., \cite[Proposition 1.5]{Langer} or
\cite[Theorem 3.1]{Trunk2}) that this implies the existence of exactly
one  eigenvalue $\lambda$ of \eqref{prob} in $\mbR$  or the upper
half-plane $\mbC^+$ and that if $\lambda\in\mbR$ with eigenfunction $\phi$
then $[A\phi,\phi]=\lambda[\phi,\phi]\le 0$.
Let $\lambda$ be such an eigenvalue with eigenfunction $\phi$.
If $\lambda$ is real, then  $-\lambda=-\overline{\lll}$ is also
an eigenvalue with the eigenfunction $\overline{\phi(-\cdot)}$
by Lemma \ref{lem21} and
$$
-\lambda[\overline{\phi(-\cdot)},\overline{\phi(-\cdot)}]=\lll[\phi,\phi]\le 0
$$
by the odd symmetry of $w$.
Thus, we get that $\lambda$ and $-\lambda$ are
two such eigenvalues, which is a contradiction.
Since $\lll\in\mbC^+$ implies $-\overline{\lambda}\in\mbC^+$,
we see that $\lambda=-\overline{\lambda}$, i.e.,
$\lambda$ is purely imaginary.
The proof of Theorem \ref{thm1} is complete.\qed
\medskip

\noindent{\bf Acknowledgment.}
The authors gratefully thank the referee for his or her helpful suggestions
which lead to the improvement of the result in and simplification of the
proof of Theorem \ref{thm1}.

\end{document}